\documentclass[12pt]{amsart}
\usepackage[colorlinks=true,citecolor=blue,linkcolor=blue]{hyperref}
\usepackage{latexsym,amsmath,amssymb}
\usepackage{accents}

\usepackage{color}

\title{Regularity for a fractional $p$-Laplace equation}

\author{Armin Schikorra}
\author{Tien-Tsan Shieh}
\author{Daniel Spector}

\address{Armin Schikorra, Mathematisches Institut, Abt. f\"ur Reine Mathematik, Albert-Ludwigs-Universit\"at, Eckerstra\ss{}e 1, 79104 Freiburg im Breisgau, Germany {\tt armin.schikorra@math.uni-freiburg.de}}
\address{Tien-Tsan Shieh, National Center for Theoretical Sciences, 2F Astronomy-Mathematics Building, National Taiwan University, No. 1, Sec. 4, Roosevelt Rd., Taipei City 106, Taiwan, R.O.C. {\tt ttshieh@ncts.ntu.edu.tw}}
\address{Daniel Spector, National Chiao Tung University, Department of Applied Mathematics, 1001 Ta Hsueh Rd, 30010 Hsinchu, Taiwan, R.O.C. {\tt dspector@math.nctu.edu.tw}}

\thanks{
A.S. is supported by DFG-grant SCHI-1257-3-1 and the DFG-Heisenberg fellowship.
D.S. is supported by the Taiwan Ministry of Science and Technology under research grants 103-2115-M-009-016-MY2 and 105-2115-M-009-004-MY2.
Part of this work was written while A.S. was visiting NCTU with support from the Taiwan Ministry of Science and Technology through the Mathematics Research Promotion Center.
}

plus 6pt minus 12pt
\newtheorem{theorem}{Theorem}

\newtheorem{proposition}[theorem]{Proposition}

\theoremstyle{definition}

\newcommand{\R}{\mathbb{R}}

\newcommand{\brac}[1]{\left (#1 \right )}

\newcommand{\barint}{
\rule[.036in]{.12in}{.009in}\kern-.16in \displaystyle\int }

\newcommand{\barcal}{\mbox{$ \rule[.036in]{.11in}{.007in}\kern-.128in\int $}}

\def\mvint_#1{\mathchoice
          {\mathop{\vrule width 6pt height 3 pt depth -2.5pt
                  \kern -8pt \intop}\nolimits_{\kern -3pt #1}}%
          {\mathop{\vrule width 5pt height 3 pt depth -2.6pt
                  \kern -6pt \intop}\nolimits_{#1}}%
          {\mathop{\vrule width 5pt height 3 pt depth -2.6pt
                  \kern -6pt \intop}\nolimits_{#1}}%
          {\mathop{\vrule width 5pt height 3 pt depth -2.6pt
                  \kern -6pt \intop}\nolimits_{#1}}}

\numberwithin{theorem}{section} \numberwithin{equation}{section}

\newcommand{\lap}{\Delta }

\newcommand{\laps}[1]{(-\lap)^{\frac{#1}{2}}}

\begin{document}
\begin{abstract}
In this note we consider regularity theory for a fractional $p$-Laplace operator which arises in the complex interpolation of the Sobolev spaces, the $H^{s,p}$-Laplacian. We obtain the natural analogue to the classical $p$-Laplacian situation, namely $C^{s+\alpha}_{loc}$-regularity for the homogeneous equation.
\end{abstract}
\maketitle

\section{Introduction and main result}
In recent years equations involving what we will call the distributional $W^{s,p}$-Laplacian, defined for test functions $\varphi$ as 
\[
(-\lap)^s_p u[\varphi] := \int_{\mathbb{R}^d} \int_{\mathbb{R}^d} \frac{|u(x)-u(y)|^{p-2}(u(x)-u(y)) (\varphi(x)-\varphi(y))}{|x-y|^{d+sp}}\;dy\, dx,
\]
have received a lot of attention, e.g. \cite{Brasco-Lindgren-2015,CastroKuusiPalatucciLocalBehaviour,CastroKuusiPalatucciJFA14,KorvenpaaKuusiPalatucci-2016,SireKuusiMingioneFracGehring,SireKuusiMingioneSelfImproving,Schikorra-CPDE}. The $W^{s,p}$-Laplacian $(-\lap)^s_p$ appears when one computes the first variation of certain energies involving the $W^{s,p}$ semi-norm
\begin{align}
[u]_{W^{s,p}(\mathbb{R}^d)}:= \brac{\int_{\mathbb{R}^d} \int_{\mathbb{R}^d} \frac{|u(x)-u(y)|^{p}}{|x-y|^{d+sp}}\;dy\, dx}^{\frac{1}{p}},\label{eq:gagliardo}
\end{align}
which was introduced by Gagliardo \cite{Gagliardo-1957} and independently by Slobodeckij \cite{Slobodecki-1958} to describe the trace spaces of Sobolev maps.

Regularity theory for equations involving this fractional $p$-Laplace operator is a very challenging open problem and only partial results are known:  
$C^{0,\alpha}_{loc}$-regularity for suitable right-hand-side data was obtained by Di Castro, Kuusi and Palatucci \cite{CastroKuusiPalatucciLocalBehaviour,CastroKuusiPalatucciJFA14}; A generalization of the Gehring lemma was proven by Kuusi, Mingione and Sire \cite{SireKuusiMingioneFracGehring,SireKuusiMingioneSelfImproving}; A stability theorem similar to the Iwaniec stability result for the $p$-Laplacian was established by the first-named author \cite{Schikorra-2015}. The current state-of-the art with respect to regularity theory is higher Sobolev-regularity by Brasco and Lindgren \cite{Brasco-Lindgren-2015}.

Aside from their origins as trace spaces, the fractional Sobolev spaces 
\[
W^{s,p}(\mathbb{R}^d):=\left\{ u \in L^p(\mathbb{R}^d) : [u]_{W^{s,p}(\mathbb{R}^d)}<+\infty\right\}
\]
also arise in the real interpolation of $L^p$ and $\dot{W}^{1,p}$.  If one alternatively considers the complex interpolation method, one is naturally led to another kind of fractional Sobolev space $H^{s,p}(\mathbb{R}^d)$, where taking the place of  the differential energy \eqref{eq:gagliardo} one can utilize the semi-norm
\begin{equation}\label{eq:hsp}
 [u]_{H^{s,p}(\mathbb{R}^d)} := \brac{\int_{\mathbb{R}^d} |D^s u|^p }^{\frac{1}{p}}.
\end{equation}
Here $D^s = (\frac{\partial^s}{\partial x^s_1},\ldots, \frac{\partial^s}{\partial x^s_d})$ is the fractional gradient for
\begin{align*}
\frac{\partial^s u}{\partial x^s_i}(x):= c_{d,s}\int_{\mathbb{R}^d} \frac{u(x)-u(y)}{|x-y|^{d+s}} \frac{x_i-y_i}{|x-y|}\;dy, \quad i = 1,\dots,d.
\end{align*}
Composition formulae for the fractional gradient have been studied in the classical work \cite{Horvath-1959}, while more recently they have been considered by a number of authors  \cite{BilerImbertKarch-2015,Caffarelli-Vazquez-2013,Caffarelli-Vazquez-2011,Schikorra-gradient,Schikorra-eps,Shieh-Spector-2014}. 
While it is common in the literature (for example in \cite{Mazya-book}) to see $H^{s,p}(\mathbb{R}^d)$ equipped with the $L^p$-norm of the fractional Laplacian $\laps{s}$ (see Section~\ref{s:proof} for a definition), we here utilize \eqref{eq:hsp} because it preserves the structural properties of the spaces for $s \in (0,1)$ more appropriately.  In particular, for $s = 1$ we have $D^1 = D$ (the constant $c_{d,s}$ tends to zero as $s$ tends to one), while for $s \in (0,1)$ the fractional Sobolev spaces defined this way support a fractional Sobolev inequality in the case $p=1$  (see \cite{SSvS-2015}).  
Let us also remark that for $p=2$ these spaces are the same, $W^{s,2} = H^{s,2}$, but for $p \neq 2$ this is not the case.

Returning to the question of a fractional $p$-Laplacian, in the context of $H^{s,p}(\mathbb{R}^d)$ computing the first variation of energies involving the $H^{s,p}$ semi-norms \eqref{eq:hsp} yields an alternative fractional version of a $p$-Laplacian, we shall call it the $H^{s,p}$-Laplacian
\[
\operatorname{div}_s (|D^s u|^{p-2}\, D^s u) = \sum_{i=1}^d \frac{\partial^s }{\partial x^s_i}(|D^s u|^{p-2}\, \frac{\partial^s u}{\partial x^s_i}).
\]

Somewhat surprisingly while the regularity theory for the homogeneous equation of the $W^{s,p}$-Laplacian
\[
 (-\lap)^s_p u = 0
\]
is far from being understood, the regularity for the $H^{s,p}$-Laplacian
\begin{equation}\label{eq:fracplapeq}
 \operatorname{div}_s (|D^s u|^{p-2}\, D^s u) = 0
\end{equation}
actually follows the classical theory, which is the main result we prove in this note:

\begin{theorem}\label{th:regularity}
Let $\Omega \subset \R^d$ be open, $p \in (2-\frac{1}{d},\infty)$ and $s \in (0,1]$. Suppose $u \in H^{s,p}(\R^d)$ is a distributional solution to \eqref{eq:fracplapeq},
that is
\begin{equation}\label{eq:fracplapeqdist}
\int_{\mathbb{R}^d} |D^su|^{p-2}D^su\cdot D^s\varphi =  0 \quad \forall \varphi \in C_c^\infty(\Omega).
\end{equation}
Then $u \in C^{s+\alpha}_{loc}(\Omega)$ for some $\alpha > 0$ only depending on $p$.
\end{theorem}

The key observation for Theorem~\ref{th:regularity} is that $v := I_{1-s}u$, where $I_{1-s}$ denotes the Riesz potential, actually solves an inhomogeneous classical $p$-Laplacian equation with good right-hand side.
\begin{proposition}\label{pr:reduction}
Let $u$ be as in Theorem~\ref{th:regularity}. Then $v := I_{1-s}u$ satisfies
\[
-\operatorname{div}(|Dv|^{p-2}Dv) \in L^\infty_{loc}(\Omega).
\]
\end{proposition}
Therefore, Theorem~\ref{th:regularity} follows from the regularity theory of the classical $p$-Laplacian: By Proposition~\ref{th:regularity}, $v$ is a distributional solution to
\[
 \operatorname{div} (|D v|^{p-2}\, D v) = \mu
\]
and $\mu$ is sufficiently integrable whence $v \in C^{1,\alpha}_{loc}(\Omega)$ \cite{DiBenedetto-1983,Evans-1982,Uralceva-1968} (see also the excellent survey paper by Mingione~\cite{Mingione-2014}).  In particular, one can apply the potential estimates by Kuusi and Mingione \cite[Theorem 1.4, Theorem 1.6]{Kuusi-Mingione-2012} to deduce that $Dv \in C^{0,\alpha}_{loc}(\Omega)$, which implies that $u \in C^{s+\alpha}_{loc}(\Omega)$.

Let us also remark, that the reduction argument used for Proposition~\ref{pr:reduction} extends the class of fractional partial differential equations introduced in~\cite{Shieh-Spector-2014}, which will be treated in a forthcoming work.

\section{Proof of Proposition~\ref{pr:reduction}}\label{s:proof}
With $\laps{\sigma}$ we denote the fractional Laplacian
\[
 \laps{\sigma} f(x):= \tilde{c}_{d,\sigma} \int_{\R^d} \frac{f(x)-f(y)}{|x-y|^{d+\sigma}}\ dy,
\]
and with $I_{\sigma}$ its inverse, the Riesz potential. Let $v := I_{1-s}u$ where $u$ satisfies \eqref{eq:fracplapeqdist}, so that 
\begin{equation}
\int_{\mathbb{R}^d} |Dv|^{p-2} Dv \cdot D^s \varphi = 0 \quad \mbox{for all $\varphi \in C^\infty_c(\Omega)$.}\label{eq:reduction}
\end{equation}
Now let $\Omega_1 \Subset \Omega$ be an arbitrary open set compactly contained in $\Omega$, and let $\phi$ be a test function supported in $\Omega_1$. Pick an open set $\Omega_2$ so that $\Omega_1\Subset \Omega_2 \Subset \Omega$ and a cutoff function $\eta$, supported in $\Omega$ and constantly one in $\Omega_2$.
Then in particular one can take \[\varphi := \eta (-\Delta)^{\frac{1-s}{2}}\phi\] as a test function in \eqref{eq:reduction} to obtain
\begin{align*}
\int_{\mathbb{R}^d} |Dv|^{p-2} Dv \cdot D^s (\eta (-\Delta)^{\frac{1-s}{2}}\phi) = 0.
\end{align*}
That is,
\begin{align*}
\int_{\mathbb{R}^d} |Dv|^{p-2} Dv \cdot D \phi = \int_{\mathbb{R}^d} |Dv|^{p-2} Dv \cdot D^s (\eta^c (-\Delta)^{\frac{1-s}{2}}\phi).
\end{align*}
where $\eta^c :=(1-\eta)$.  We set 
\[
 T(\phi) := D^s (\eta^c (-\Delta)^{\frac{1-s}{2}}\phi) 
\]

Now we show that by the disjoint support of $\eta^c$ and $\phi$ we have 
 \begin{equation}\label{eq:finalgoal}
  \|T(\phi)\|_{L^p(\R^d)} \leq C_{\Omega_1,\Omega_2,d,s,p}\, \|\phi\|_{L^1(\R^d)}. 
 \end{equation}
Once we have this, the claim is proven as H\"older's inequality and realizing the $L^\infty$ norm via duality implies 
 \[
 -\operatorname{div}(|Dv|^{p-2}Dv) \in L^\infty_{loc}(\Omega).
 \]
To see \eqref{eq:finalgoal}, we use the disjoint support arguments as in \cite[Lemma A.1]{Blatt-Reiter-Schikorra-2016} \cite[Lemma 3.6.]{Martinazzi-Schikorra-Maalaoui-2015}: First we see that since $\eta^c(x) \phi(x) \equiv 0$,
\[
 T(\phi) = \tilde{c}_{d,1-s}\, D^s \int_{\mathbb{R}^d}  \frac{ \eta^c(x)\phi(y)}{|x-y|^{N+1-s}}\;dy.
\]
Now taking a cutoff-function $\zeta$ supported in $\Omega_2$, $\zeta \equiv 1$ on $\Omega_1$ we have
\[
  T(\phi)= \tilde{c}_{d,1-s}\, D^s \int_{\mathbb{R}^d}  \frac{ \eta^c(x)\zeta(y) \phi(y)}{|x-y|^{N+1-s}}\;dy =  \tilde{c}_{d,1-s}\,\int_{\mathbb{R}^d}  k(x,y)\,  \phi(y)\;dy,
\]
where
\[
 k(x,y) := D^s_x \kappa(x,y) := D^s_x \frac{ \eta^c(x)\, \zeta(y) }{|x-y|^{N+1-s}}.
\]
The positive distance between the supports of $\eta^c$ and $\zeta$ implies that these kernels $k$, $\kappa$ are a smooth, bounded, integrable (both, in $x$ and in $y$), and thus by a Young-type convolution argument we obtain \eqref{eq:finalgoal}. One can also argue by interpolation,
\begin{align*}
\left\| \int_{\mathbb{R}^d}  \kappa(x,y)\, \phi(y) \right\|_{L^p(\mathbb{R}^d)} \leq \|\phi\|_{L^1(\mathbb{R}^d)},
\end{align*}
as well as
\begin{align*}
\left\| \int_{\mathbb{R}^d}  D_x \kappa(x,y)\, \phi(y)\, dy \right\|_{L^p(\mathbb{R}^d)} \leq \|\phi\|_{L^1(\mathbb{R}^d)}.
\end{align*}
Interpolating this implies the desired result that
\begin{align*}
\left\| \int_{\mathbb{R}^d}  D^s_x \kappa(x,y)\, \phi(y)\, dy \right\|_{L^p(\mathbb{R}^d)} \leq \|\phi\|_{L^1(\mathbb{R}^d)}.
\end{align*}
Thus \eqref{eq:finalgoal} is established and the proof of Proposition~\ref{pr:reduction} is finished.
\qed

\bibliographystyle{amsplain}
\bibliography{bib}

\end{document}